\documentclass[12pt]{amsart}      %{article} was 12pt latex e
\usepackage{amssymb}
\usepackage{eucal}
\usepackage{graphicx}
\usepackage{amsmath}
\usepackage{amscd}
\usepackage[all]{xy}           %xypic macro for latex2.09

\usepackage{amsfonts,latexsym}
\newtheorem{theorem}{Theorem}[section]
\newtheorem{prop}[theorem]{Proposition}
\newtheorem{defn}[theorem]{Definition}
\newtheorem{lemma}[theorem]{Lemma}
\newtheorem{coro}[theorem]{Corollary}
\newtheorem{prop-def}{Proposition-Definition}[section]

\newcommand{\nc}{\newcommand}
\renewcommand{\Bbb}{\mathbb}
\renewcommand{\frak}{\mathfrak}
\newcommand{\efootnote}[1]{}
\renewcommand{\textbf}[1]{}
\newcommand{\delete}[1]{}
\nc{\dfootnote}[1]{{}}          %{{}}
\nc{\ffootnote}[1]{\dfootnote{#1}}
\nc{\mfootnote}[1]{\footnote{#1}} % Use this to show footnotes
\nc{\ofootnote}[1]{\footnote{\tiny Older version: #1}} % Use this to show footnotes

\nc{\mlabel}[1]{\label{#1}}  % Use this to suppress names
\nc{\mcite}[1]{\cite{#1}}  % Use this to suppress names
\nc{\mref}[1]{\ref{#1}}  % Use this to suppress names
\nc{\mbibitem}[1]{\bibitem{#1}} % Use this to show number name

\nc{\bin}[2]{ (_{\stackrel{\scs{#1}}{\scs{#2}}})}  %binomial coeff
\nc{\binc}[2]{ \left (\!\! \begin{array}{c} \scs{#1}\\
    \scs{#2} \end{array}\!\! \right )}  %binomial coeff
\nc{\bincc}[2]{  \left ( {\scs{#1} \atop
    \vspace{-1cm}\scs{#2}} \right )}  %binomial coeff
\nc{\bs}{\bar{S}}
\nc{\cosum}{\sqsubset}
\nc{\la}{\longrightarrow}
\nc{\rar}{\rightarrow}
\nc{\dar}{\downarrow}
\nc{\dap}[1]{\downarrow \rlap{$\scriptstyle{#1}$}}
\nc{\uap}[1]{\uparrow \rlap{$\scriptstyle{#1}$}}
\nc{\defeq}{\stackrel{\rm def}{=}}
\nc{\disp}[1]{\displaystyle{#1}}
\nc{\dotcup}{\ \displaystyle{\bigcup^\bullet}\ }
\nc{\hcm}{\ \hat{,}\ }
\nc{\hts}{\hat{\otimes}}
\nc{\free}[1]{\bar{#1}}
\nc{\hcirc}{\hat{\circ}}
\nc{\lleft}{[}
\nc{\lright}{]}
\nc{\curlyl}{\left \{ \begin{array}{c} {} \\ {} \end{array}
    \right .  \!\!\!\!\!\!\!}
\nc{\curlyr}{ \!\!\!\!\!\!\!
    \left . \begin{array}{c} {} \\ {} \end{array}
    \right \} }
\nc{\longmid}{\left | \begin{array}{c} {} \\ {} \end{array}
    \right . \!\!\!\!\!\!\!}
\nc{\ora}[1]{\stackrel{#1}{\rar}}
\nc{\ola}[1]{\stackrel{#1}{\la}}%${\Bbb Z}$
\nc{\ot}{\otimes}
\nc{\sprod}{\bullet}
\nc{\scs}[1]{\scriptstyle{#1}}
\nc{\mrm}[1]{{\rm #1}}
\nc{\margin}[1]{\marginpar{\rm #1}}   %{\rm #1}}
\nc{\dirlim}{\displaystyle{\lim_{\longrightarrow}}\,}
\nc{\invlim}{\displaystyle{\lim_{\longleftarrow}}\,}
\nc{\mvp}{\vspace{0.3cm}}
\nc{\tk}{^{(k)}}
\nc{\tp}{^\prime}
\nc{\ttp}{^{\prime\prime}}
\nc{\svp}{\vspace{2cm}}
\nc{\vp}{\vspace{8cm}}
\nc{\proofbegin}{\noindent{\bf Proof: }}
\nc{\proofend}{$\blacksquare$ \vspace{0.3cm}}
\nc{\modg}[1]{\!<\!\!{#1}\!\!>}
\nc{\intg}[1]{F_C(#1)}
\nc{\lmodg}{\!<\!\!}
\nc{\rmodg}{\!\!>\!}
\nc{\cpi}{\widehat{\Pi}}
\nc{\sha}{{\mbox{\cyr X}}}  %used to be \cyr
\nc{\shap}{{\mbox{\cyrs X}}} %sha as product
\nc{\shpr}{\diamond}    %Shuffle product
\nc{\shplus}{\shpr^+}
\nc{\shprc}{\shpr_c}    %Cartier's product
\nc{\vep}{\varepsilon}
\nc{\labs}{\mid\!}
\nc{\rabs}{\!\mid}
\nc{\rchar}{\mrm{char}}
\nc{\Fil}{\mrm{Fil}}
\nc{\Hom}{\mrm{Hom}}
\nc{\id}{\mrm{id}}
\nc{\im}{\mrm{im}}
\nc{\incl}{\mrm{incl}}
\nc{\map}{\mrm{Map}}
\nc{\mchar}{\rm char}
\nc{\supp}{\mathrm Supp}

\nc{\Alg}{\mathbf{Alg}}
\nc{\Bax}{\mathbf{Bax}}
\nc{\bff}{\mathbf f}
\nc{\bfk}{{\bf k}}
\nc{\bfone}{{\bf 1}}
\nc{\bfx}{\mathbf x}
\nc{\bfy}{\mathbf y}
\nc{\base}[1]{\bfone^{\otimes ({#1}+1)}} %{{a_{#1}}}
\nc{\Cat}{\mathbf{Cat}}

\nc{\detail}{\marginpar{\bf More detail}
    \noindent{\bf Need more detail!}
    \svp}
\nc{\Int}{\mathbf{Int}}
\nc{\Mon}{\mathbf{Mon}}
\nc{\remarks}{\noindent{\bf Remarks: }}
\nc{\Rings}{\mathbf{Rings}}
\nc{\Sets}{\mathbf{Sets}}

\nc{\BA}{{\Bbb A}} \nc{\CC}{{\Bbb C}} \nc{\DD}{{\Bbb D}}
\nc{\EE}{{\Bbb E}} \nc{\FF}{{\Bbb F}} \nc{\GG}{{\Bbb G}}
\nc{\HH}{{\Bbb H}} \nc{\LL}{{\Bbb L}} \nc{\NN}{{\Bbb N}}
\nc{\KK}{{\Bbb K}} \nc{\QQ}{{\Bbb Q}} \nc{\RR}{{\Bbb R}}
\nc{\TT}{{\Bbb T}} \nc{\VV}{{\Bbb V}} \nc{\ZZ}{{\Bbb Z}}

\nc{\cala}{{\mathcal A}} \nc{\calc}{{\mathcal C}}
\nc{\cald}{{\mathcal D}} \nc{\cale}{{\mathcal E}}
\nc{\calf}{{\mathcal F}} \nc{\calg}{{\mathcal G}}
\nc{\calh}{{\mathcal H}} \nc{\cali}{{\mathcal I}}
\nc{\call}{{\mathcal L}} \nc{\calm}{{\mathcal M}}
\nc{\caln}{{\mathcal N}} \nc{\calo}{{\mathcal O}}
\nc{\calp}{{\mathcal P}} \nc{\calr}{{\mathcal R}}
\nc{\calt}{{\mathcal T}} \nc{\calw}{{\mathcal W}}
\nc{\calk}{{\mathcal K}} \nc{\calx}{{\mathcal X}}
\nc{\CA}{\mathcal{A}}

\nc{\fraka}{{\frak a}}
\nc{\frakA}{{\frak A}}
\nc{\frakB}{{\frak B}}
\nc{\frakH}{{\frak H}}
\nc{\frakM}{{\frak M}}
\nc{\frakm}{{\frak m}}
\nc{\frakP}{{\frak P}}
\nc{\frakp}{{\frak p}}

\font\cyr=wncyr10
\font\cyrs=wncyr7
\nc{\redt}[1]{\textcolor{red}{#1}}
\begin{document}

\title[Multiple zeta values and Rota--Baxter algebras]
{Multiple zeta values and Rota--Baxter algebras\\
\ \\
\normalsize\rm (Dedicated to Professor Melvyn Nathanson for his 60th birthday)}%
%=========================================================================================
%
\author{Kurusch Ebrahimi-Fard}

\address{I.H.\'E.S.
         Le Bois-Marie,
         35, Route de Chartres,
         F-91440 Bures-sur-Yvette,
         France}
\email{kurusch@ihes.fr}
\author{Li Guo}

\address{Department of Mathematics and Computer Science,
         Rutgers University,
         Newark, NJ 07102}
\email{liguo@newark.rutgers.edu}

%=========================================================================================

%\date{\today}
%=========================================================================================
%=========================================================================================
%=========================================================================================

%\begin{document}
\maketitle

%=========================================================================================
\begin{abstract}
We study multiple zeta values and their generalizations from the point of
view of Rota--Baxter algebras. We obtain a general framework for this
purpose and derive relations on multiple zeta values from
relations in Rota--Baxter algebras.
\end{abstract}
%=========================================================================================
%=========================================================================================

%\tableofcontents
%=========================================================================================

\setcounter{section}{0}

%=========================================================================================
\section{Introduction}

The purpose of this paper is to establish the relationship between
Rota--Baxter algebras and multiple zeta values, multiple polylogarithms
and their $q$-analogs.
\medskip

%\subsection{History of MZVs}
Multiple zeta values, henceforth abbreviated MZVs, are defined by
\begin{equation}
  \zeta(s_1,\cdots,s_k):=\sum_{n_1>n_2>\cdots>n_k\geq 1} \frac{1}{n_1^{s_1}\cdots n_k^{s_k}}
  \mlabel{MZV}
\end{equation}
where the $s_i$ are positive integers with $s_1>1$.
%{We call $k$ the length of $\zeta(s_1,\cdots,s_k)$ and
%$\sum_{i=1}^{k} s_i$ its depth.}

The earliest algebraic relation among multiple zeta values
traces back to Euler. Their systematic study started in early
1990s with the work of Hoffman~\mcite{Ho0} and Zagier~\mcite{Za}.
Since then MZVs and their generalizations have been studied
intensively by numerous authors with
connections to arithmetic geometry, mathematical physics, quantum
groups and knot theory. Surveys on the related work can be found
in~\mcite{B-B-B-L,Ho3,Ca1,Wa1,Wa2,Zud}.
Lately generalizations of MZVs, such as multiple polylogarithms (MPLs)
have also been shown to be important, in both pure mathematics
\mcite{Za,Ca1} and theoretical physics \mcite{B-K,Krei,M-U-W}.
Investigations of their possible q-analogs, especially with
respect to algebraic aspects were started in \mcite{Brad,Sch,Zh,Zud}.
\medskip

%=========================================================================================
%\subsection{History of Rota--Baxter algebras and mixable shuffles}

A Rota--Baxter operator of weight $\lambda$ is a linear operator
$P$ on an algebra $A$ such that
\begin{equation}
  P(x)P(y)=P(xP(y))+P(P(x)y) + \lambda P(xy),\ x,y \in A.
  \mlabel{RB0}
\end{equation}
Here $\lambda$ is a fixed constant in the base ring. \delete{The
image of $P$ forms a subalgebra in $A$.}
Rota--Baxter algebra was first introduced by Baxter~\cite{Ba} in
1960 to study the theory of fluctuations in probability. It was
further studied in the next two decades by a number of
mathematicians, especially Rota who greatly contributed to the
study of the Rota--Baxter algebra by his pioneer{ing} work in
the late 1960s and early 1970s~\cite{Ro1, Ro2, Ro3} and by his
survey articles in late 1990s~\cite{Ro4,Ro5}.

In the last few years there have been further developments in
Rota-Baxter algebras with
applications to quantum field theory~\cite{C-K1, C-K2, E-G-K1,
E-G-K2,E-G-K3,E-G10,E-G-G-V}, dendriform algebras~\cite{Ag1,Ag2,A-L,EF,
E-G2,E-G4}, number theory~\cite{Gu3}, Hopf algebras~\cite{A-G-K-O,E-G1} and
combinatorics~\cite{Gu2}. Key to some of these developments in
Rota--Baxter algebra is the realization of the free objects in
which the product is defined by mixable
shuffles~\cite{G-K1,G-K2}.
\medskip
%=========================================================================================
%\subsection{Motivation}

All known algebraic relations among MZVs are given by
the combination of the shuffle product of the integral representation
of MZVs and the stuffle (i.e., quasi-shuffle) product of the sum
representation of the MZVs, and their degenerated forms.
Conjecturally, all algebraic relations among MZVs can be
obtained this way~\mcite{I-K-Z}.
These products are also the shuffle product and mixable shuffle product
in the free commutative Rota--Baxter algebras~\mcite{G-K1,G-K2,E-G1},
as we will elaborate further below.
It is therefore reasonable to
expect that much of the recent work on algebraic relations for
MZVs can be viewed and expanded in the framework of Rota--Baxter
algebras. The current paper is a first step in this direction.
\medskip

In order to make precise the connection between Rota--Baxter
algebras and MZVs, and {to} set up a general framework to
deal with the various generalizations of MZVs, we review the
constructions of free commutative Rota--Baxter algebras in
Section~\mref{sec:free} and their relations with shuffle type
products {for} MZVs.
%we view MZVs and their generalizations
%in the framework of Rota--Baxter algebra which allows us to give a uniform
%approach and generalizations to these variations.
%In Section~\ref{sec:free}
%we describe the equivalence between shuffle type products in Rota--Baxter
%algebras and similar products in MZVs,
These include the relation between mixable shuffles in
Rota--Baxter algebras and quasi-shuffles and generalized shuffles
in MZVs, and the relation between the product of Cartier and stuffle
product of MZVs. In Section~\mref{sec:setup}, we use the language of free
Rota--Baxter algebras to define the concept of MZV algebras which
will include as special cases the MZVs, MPLs and $q$-MZVs. In
Section~\ref{sec:app}, this setup is used
to derive identities in MZVs from identities in Rota--Baxter
algebras.

%\subsection{Future work}
%We also use Rota--Baxter operators as a guide to define $q$-MZVs
%and to study their properties. In a future work, we will study
%Rota--Baxter operators in the context of motivic theories for
%MZVs~\cite{Ca1,Te}.

%=========================================================================================
\section{Free Rota-Baxter algebras and double shuffle}
\mlabel{sec:free} We start with reviewing the concept of
free Rota--Baxter algebra because it most precisely and broadly
relates Rota-Baxter algebras {to} MZVs. On one hand, the products in free
Rota--Baxter algebras turn out to be the same as the products for
MZVs. On the other hand, free Rota--Baxter algebras provide a
general framework to define and study MZVs and their
generalizations. Eventually, various classes of MZVs will be shown
to be subquotients of free commutative Rota--Baxter algebras.

\subsection{Rota--Baxter algebras}
We first fix some notations. Throughout this paper, we will only
consider commutative rings and algebras. For noncommutative
Rota--Baxter algebras and applications to physics and operads,
see~\mcite{E-G4,E-G10,E-G11,E-G-K3}. Let $\bfk$ be a unitary ring, that is,
a ring with an identity which we denote by 1,
and let $\lambda\in \bfk$ be fixed. A unitary (resp. nonunitary) Rota--Baxter
$\bfk$-algebra {(RBA)} of weight $\lambda$ is a pair $(R,P)$
in which $R$ is a unitary (resp. nonunitary) $\bfk$-algebra and
$P: R \to R$ is a $\bfk$-linear map such that
\begin{equation}
 P(x)P(y) = P(xP(y))+P(P(x)y)+ \lambda P(xy),\ \forall x,\ y\in R.
\mlabel{eq:Ba}
\end{equation}

We will focus on two Rota-Baxter operators here. One is the
integration operator
\begin{equation}
 I(f)(x)=\int_0^x f(t)dt
 \mlabel{eq:int}
 \end{equation}
defined on continuous functions $f(x)$ on $[0,\infty)$.
Then the integration by parts formula reads
$$ I(fI(g))=I(f)I(g)-I(I(f)g)$$
showing that the integration operator is a Rota-Baxter operator of
weight 0. The second operator is the sum operator that we will discuss
at the beginning of the next section.

Let $A$ be a unitary $\bfk$-algebra. A unitary
Rota--Baxter algebra $(F(A),P_A)$ of weight $\lambda$ is called a
free unitary Rota--Baxter algebra over $A$ if there is a unitary algebra
homomorphism $j_A:A\to F(A)$ with the property that, for any
unitary Rota--Baxter $\bfk$-algebra $(R,P)$ of weight $\lambda$ and
unitary algebra homomorphism $f: A\to R$, there is a unitary Rota--Baxter
$\bfk$-algebra homomorphism $\free{f}: (F(A),P_A)\to (R,P)$ such
that $f=\free{f} \circ j_A$, in other words, such that the diagram
\[\xymatrix{
A \ar[rr]^(0.4){j_A} \ar[drr]_{f}
    && F(A) \ar[d]^{\free{f}} \\
&& R } \]
commutes.
When all unitary is replaced by nonunitary in the above definition, we
obtain the concept of the free nonunitary Rota-Baxter algebra over $A$.

\subsection{Free Rota-Baxter algebras and mixable shuffle product}

For a given unitary commutative algebra $A$, define
$\sha^+(A):=\bigoplus_{n\geq 1} A^{\otimes n}$.
We briefly recall the definition of {\bf mixable shuffle product $\shplus$}
on $\sha^+(A)$. For details, see~\mcite{E-G1,G-K1,G-K2,Gu1}.

Consider two pure tensors $a:=a_1\otimes \ldots \otimes a_m\in
A^{\ot m}$ and $b:=b_1\otimes \ldots \otimes b_n \in A^{\ot n}$.
As is well-known, a {\bf shuffle} of $a$ and $b$ is a tensor
permutation of $a_i$ and $b_j$ without chang{ing} the order
of the $a_i$s and $b_j$s, and the shuffle product $a \shap b$ of $a$
and $b$ is the sum of shuffles of $a$ and $b$. For example
$$
    a_1 \shap (b_1\otimes b_2)
        = a_1\otimes b_1\otimes b_2 + b_1\otimes a_1\otimes b_2
            + b_1\otimes b_2\otimes a_1.
$$
More generally, a {\bf mixable shuffle} is a shuffle in which some
pairs $a_i\otimes b_j$ (but not $b_j\ot a_i$) are merged into
$\lambda a_i b_j$. The mixable shuffle product $a \shplus b$ of $a$
and $b$ is the sum of mixable shuffles of $a$ and $b$.
For example, \allowdisplaybreaks{
\begin{eqnarray}
    a_1 \shplus (b_1\otimes b_2)
            &=& a_1\otimes b_1\otimes b_2 + b_1\otimes a_1\otimes b_2
                        + b_1\otimes b_2\otimes a_1\ \  ({\rm shuffles}) \\
            & & + \lambda a_1b_1\otimes b_2 +\lambda b_1\otimes a_1b_2 \ \
                        ({\rm merged\ shuffles}). \
 \mlabel{eq:mixShuf}
\end{eqnarray}}
So the {mixable} shuffle product is the shuffle product when
$\lambda=0$. With the product $\shplus$, $\sha^+(A)$ is a
nonunitary commutative algebra. Let $\bfk\oplus \sha^+(A)$
be the unitary algebra after unitarization and let
$$\sha(A):= A\otimes (\bfk\oplus \sha^+(A))$$
be the tensor product algebra with its product denoted by $\shpr$.
So for $a\ot (u+a')$ and $b\ot (v+b')$ in $\sha(A)$ with
$a,b\in A$, $u,v\in \bfk$ and $a',b'\in \sha^+(A)$, we have the product
\begin{equation}
 (a\ot a')\shpr (b\ot b')=(ab) \ot ( uv+va'+ub'+ a'\shplus b').
 \mlabel{eq:shpr}
\end{equation}

We have the natural unitary algebra homomorphism $j_A:A\to
\sha(A)$ sending $a\in A$ to $a\ot 1\in \sha(A)$. The following
theorem is proved in~\mcite{G-K1,G-K2}.

\begin{theorem} The algebra $\sha(A)$ with the shift operator
$P_A:\sha(A)\to \sha(A), P_A(a):=1\otimes a$ is the free
commutative unitary Rota--Baxter algebra over $A$. When $A=\bfk[X]$, it is
the free commutative unitary Rota-Baxter algebra over $X$.
\end{theorem}

If $A$ is a nonunitary algebra, then the free commutative nonunitary algebra
over $A$ can be constructed as a subalgebra of $\sha(\tilde{A})$~\mcite{G-K2}.
Here $\tilde{A}=\bfk\oplus A$ is the unitarization of $A$.

\subsection{Connection with quasi-shuffle}
It is shown in~\mcite{E-G1} that the mixable shuffle product can be
recursively defined as follows.
For any $m,n\geq 1$ and
$a=a_1\ot \cdots \ot a_m\in A^{\ot m}$,
$b=b_1\ot \cdots \ot b_n\in A^{\ot n}$,
then
{\small
\begin{equation}
 a\shplus b = \left \{ \begin{array}{ll} a_1\ot b_1 + b_1\ot a_1 + \lambda a_1b_1,& \ m=n=1,
 \\[0.2cm]
 a_1\ot  b_1\ot  \cdots\ot  b_n + b_1\ot \big(a_1\shplus (b_2\ot \cdots\ot b_n)\big) \\
  \hspace{5cm} + \lambda (a_1 b_1)\ot  b_2\ot \cdots\ot  b_n, & m=1, n\geq 2,
  \\[0.2cm]
 a_1\ot \big ((a_2\ot \cdots\ot  a_m)\shplus b_1 \big) + b_1\ot a_1\ot \cdots\ot a_m  \\
  \hspace{6cm} + \lambda (a_1b_1)\ot  a_2\ot \cdots\ot  a_m,  & m\geq 2, n=1,
  \\[0.2cm]
 a_1\ot  \big ((a_2\ot \cdots\ot a_m)\shplus (b_1\ot \cdots\ot  b_n)\big )\\
  \hspace{2cm} + b_1\ot  \big ((a_1\ot  \cdots\ot  a_m)\shplus (b_2\ot \cdots\ot  b_n)\big) & \\
  \hspace{2.5cm}  + \lambda (a_1 b_1)\ot  \big ( (a_2\ot \cdots\ot a_m) \shplus
    (b_2\ot \cdots\ot  b_n)\big ), &  m, n\geq 2.
\end{array}
\right .
\mlabel{eq:quasi2}
\end{equation}
}

This is a mild generalization of the quasi-shuffle product introduced by
Hoffman~\mcite{Ho2} at about the same time as the mixable shuffle product to study MZVs.
We briefly recall his construction.
Let $X$ be a set with a grading given by finite subsets $X_n, n\geq 1,$
(locally finite set)
and with either a graded associative commutative product $[\cdot, \cdot]: X\times X\to X$
or the zero product $[\cdots,\cdots]:X\times X \to \{0\}.$
Let $\frakA$ be the noncommutative free algebra on $X$ with identity 1.
Elements in $X$ are called letters and monomials in $\frakA$ are
called words.
Hoffman's {\bf quasi-shuffle product} is
the product $*$ on $\frakA$ recursively defined by
\begin{enumerate}
\item  $1*w=w*1=w$ for any word $w$;
\item $(aw_1) * (bw_2)=a(w_1*(bw_2))+b((aw_1)*w_2)+[a,b](w_1*w_2)$,
       for any words $w_1,w_2$ and letters $a,b$.
\end{enumerate}
Then $(\frakA,*)$ is a commutative algebra.
We thus have~\mcite{E-G1}
\begin{theorem} The products $*$ coincides with
$\shplus$ on $\sha^+(A)$ in the special case when $\lambda=1$ and when
$A$ is the algebra $\bfk\oplus (\oplus_{x\in X} \bfk x)$ with product given
by $[\cdot,\cdot]$.
\end{theorem}

The mixable shuffle product was also found by Goncharov~\mcite{Go3}
in the context of motivic shuffle relations.
\smallskip

As another interesting link between Rota--Baxter algebras and
MZVs, we note that a construction of free Rota--Baxter
algebras over a set was obtained over 30 years ago by
Cartier~\mcite{Ca2} where the product is defined in terms of
ordered subsets. Recently~\mcite{Brad} the stuffle product for
MZVs was described in a similar way using order preserving
injections instead of ordered subsets.

\section{Rota--Baxter algebra setup of MZVs and their generalizations}
\mlabel{sec:setup}

The above section gives strong evidence on the intrinsic
connection between Rota--Baxter algebra and MZVs. In order to effectively
apply results of Rota--Baxter algebras to MZVs, it is desirable to
give a Rota--Baxter algebra structure on the set of MZVs.
As we see in Eq (\mref{MZV}), the MZVs are defined by iterated sums, given
by iterations of the sum operator
$$P(f)(x):= \sum_{n\geq 1} f(x+n).$$
Under certain convergency conditions, such as $f(x) = O(x^{-2})$
and $g(x)=O(x^{-2})$,
$P(f)(x)$ and $P(g)(x)$ are defined by absolutely convergent
series and we have
\begin{eqnarray}
\lefteqn{P(f)(x) P(g)(x) = \sum_{m\geq 1} f(x+m) \sum_{n\geq 1} g(x+n)}\notag \\
&=& \sum_{n>m\geq 1} f(x+m) g(x+n) + \sum_{m>n\geq 1} f(x+m)g(x+n)
+ \sum_{m\geq 1} f(x+m)g(x+m) \mlabel{eq:sum1}\\
&=& P(f P(g))(x) + P(gP(f))(x) + P(fg)(x)\notag
\end{eqnarray}
since
\begin{eqnarray*}
P(fP(g))(x)&=& \sum_{m=1}^\infty f(x+m) P(g)(x+m) \\
&=& \sum_{m=1}^\infty f(x+m) \big (\sum_{k=1}^\infty g(x+m+k) \big )\\
&=& \sum_{n>m\geq 1} f(x+m) g(x+n).
\end{eqnarray*}

This shows that the operator $P$ is a Rota-Baxter operator of weight 1
on certain functions. However, the sum operator and its
iteration are not defined on some other functions. So
we can only expect that subsets of a Rota--Baxter algebra can be
applied to the MZV study. This motivates the following construction.

\subsection{MZV algebras}

Let $R$ be a $\bfk$-algebra. Let $P$ be a partially defined map
from $R$ to $R$. We call $P$ a partially defined Rota--Baxter operator if
\begin{equation}
P(f)P(g)=P(fP(g))+P(P(f)g)+\lambda P(fg)
\mlabel{eq:RB3}
\end{equation}
if all terms are defined.
%For example,
%let $R$ be the $\RR$-algebra of continuous functions on
%$(0,\infty)$ and formally define
%$$P(f)(x)=\sum_{n\geq 1} f(x+n).$$
%We called the sum {\bf well-defined} if it converges absolutely.
%Then $P$ is a Rota--Baxter operator in this more general sense.

For $\bff:=(f_1,\cdots,f_n)\in R^n$, formally define
$$ P_\bff:=P(f_1P(f_2\cdots P(f_n)\cdots )).$$
%To abbreviate the notation $(f_1,\cdots,f_n)$, we use
%$\disp{\cosum_{i=1}^n x_i}$ in place of the notation $\disp{\Cat_{i=1}^n x_i}$
%used in~\cite{B-B-B-L}.
We define a {\bf filtered $\bfk$-algebra} to be a nonunitary $\bfk$-algebra
$A$ with a decreasing sequence
$A_n, n\geq 0$, of ideals such that $A_mA_n\subseteq A_{m+n}$.

\begin{defn}
Let $R$ be a $\bfk$-algebra with a partially defined Rota--Baxter
operator $P$. A filtered subalgebra $A$ of $R$ is called {\bf
iteratedly summable} (of level $k$) if the formal symbols
$P_\bff$ are well-defined for all $\bff\in A^n$, $n\geq 1$, with
$f_n\in A_k$.
For an iteratedly summable subalgebra $A$ (of level $k$), we call the set
$$
 \frakA_k:=\{ P_{(f_1,\cdots,f_n)} \big | f_i\in A, 1\leq i\leq n, f_n\in A_k\}
$$
the {\bf MZV algebra} generated by $A$.
\end{defn}
We will show below (Theorem~\mref{thm:MZVsubAlg})
that $\frakA_k$ is indeed an algebra.

\subsection{The abstract MZV algebra}

%We consider a nonunitary $\bfk$-algebra $A$ with a decreasing sequence
%$A_n, n\geq 0$, of ideals such that $A_mA_n\subseteq A_{m+n}$.
%We assume that $A$ has a basis $X$ and $X_n$ be a basis of $A_n$.
Let $A$ be a filtered $\bfk$-algebra. We can construct a MZV-algebra
generated by $A$ as follows.
%For example, let $X$ be the locally finite graded set in Hoffman's
%quasi-shuffle product. Then $A=\oplus_{x\in X} \bfk x$ is such an algebra.
%that $X_mX_n\subseteq X_{m+n}$. We consider more general case in order
%to cover generalizations of MZVs such as $q$-MZVs.

Let $\tilde{A}$ be the unitarization of $A$. So
$\tilde{A}=\bfk \oplus A$ with componentwise addition and with product
defined by $(m,a)(n,b)=(mn, mb+na+ab).$
In the free Rota--Baxter algebra $\sha(\tilde{A})$ with Rota-Baxter
operator $P_A$, we have $P_A(f)=1\ot f$. So
$P_{A,{(f_1,\cdots,f_n)}}=1\ot f_1\ot \cdots f_n$.
Therefore, for $k\geq 1$, the MZV algebra generated by $A$ in
$\sha(\tilde{A})$ is
the subspace $\frakM(A)_k$ of $\sha(\tilde{A})$
generated by pure tensors of the form
$1\ot a_1 \ot \cdots \ot a_n\in  1\ot A^{\ot n}$ with $a_n\in A_k$.

\delete{
\begin{proof}
We just need to prove that $\frakM(A)_k$ is closed under the mixable shuffle
product. Let
$1\ot a_1\ot \cdots \ot a_m$ and $1\ot b_1\ot \cdots \ot b_n$ be two
pure tensors in $\frakM(A)_k$. So by Eq (\mref{eq:shpr}) we have
$a_i, b_j\in A$ for $1\leq i\leq m, 1\leq j\leq n$,
and $a_m,b_n\in A_k$. Then
$$(1\ot a_1\ot \cdots \ot a_m) \shpr (1\ot b_1\ot \cdots \ot b_n)
= 1\ot \big( (a_1\ot \cdots \ot a_m)\shplus (b_1\ot \cdots \ot b_n)\big).$$
By the definition of the mixable shuffle, each mixable shuffle
of $a_1\ot \cdots \ot a_m$ and $b_1\ot \cdots b_n$ is of the form
$c_1\ot \cdots \ot c_\ell$ where $c_\ell$ is either $a_m$ or $b_n$ or
$\lambda a_m b_n$. So in either case, $c_\ell$ is in $A_k$.
Thus $(1\ot a_1\ot \cdots \ot a_m) \shpr (1\ot b_1\ot \cdots \ot b_n)$
is in $\frakM_k$.
\end{proof}
}

\begin{theorem}
Let $A$ be an iteratedly summable subalgebra of $R$ of level $k$.
\begin{enumerate}
\item
$\frakA_k$ is a subalgebra of $R$.
\mlabel{it:ak}
\item
$\frakM(A)_k$ is a subalgebra of $\sha(\tilde{A})$.
\mlabel{it:mk}
\item
There is an algebra surjection
$$
    \frakP_k: \frakM(A)_k \to \frakA_k
$$
sending $1\ot f_1\ot \cdots \ot f_n$ to $P_{(f_1,\cdots,f_n)}$.
\mlabel{it:surj}
\item
Given an evaluation, that is, an algebra homomorphism $\nu: A\to \bfk$,
we obtain an algebra homomorphism
$\nu \circ \frakP_k: \frakM(A)_k \to \bfk.$
\mlabel{it:eval}
\end{enumerate}
\mlabel{thm:MZVsubAlg}
\end{theorem}
\begin{proof}
(\mref{it:ak})
We only need to prove that $\frakA$ is closed under multiplication.
For this we just need to prove that for any
$(f_1,\cdots,f_m)\in A^m$, $(g_1,\cdots,g_n)\in A^n$ with
$f_i,g_j\in A$ and $f_m,g_n\in A_k$, the product
$P_{(f_1,\cdots,f_m)}P_{(g_1,\cdots,g_n)}$ are still of this form.
We prove this by induction on $m+n$.
When $m=n=1$, we have $f_1,g_1\in A_k$. Then
\begin{eqnarray*}
P_{f_1}P_{g_1}&=& P(f_1)P(g_1)\\
&=& P(f_1P(g_1))+P(g_1P(f_1))+\lambda P(f_1g_1).
\end{eqnarray*}
So we are done. Assuming the claim is true for $m+n\leq k$ and consider
the case of $m+n=k+1$. Then the Rota-Baxter relation~(\mref{eq:RB3})
\begin{eqnarray*}
P_{(f_1,\cdots,f_m)}P_{(g_1,\cdots,g_n)}&=&
    P(f_1 P_{(f_2,\cdots,f_m)})P(g_1P_{(g_2,\cdots,g_n)})\\
    &=& P(f_1 P_{(f_2,\cdots,f_m)}P(g_1P_{(g_2,\cdots,g_n)}))
    +P(P(f_1 P_{(f_2,\cdots,f_m)})g_1P_{(g_2,\cdots,g_n)}) \\
&&    +\lambda P(f_1 P_{(f_2,\cdots,f_m)}g_1P_{(g_2,\cdots,g_n)}).
    \end{eqnarray*}
By the induction hypothesis,
$P_{(f_2,\cdots,f_m)}P(g_1P_{(g_2,\cdots,g_n)})=
P_{(f_2,\cdots,f_m)}P_{(g_1,\cdots,g_n)}$ is a sum of terms of the form
$P_{(h_1,\cdots,h_\ell)}$ with $h_i\in A$ and $h_\ell\in A_k$.
So $P(f_1 P_{(f_2,\cdots,f_m)}P(g_1P_{(g_2,\cdots,g_n)}))$
is also a sum of the form
$P_{(f_1,h_1,\cdots,h_\ell)}$, so are still in $\frakA_k$.
The same argument applies to the other two terms. This completes
the induction.

(\mref{it:mk}) By construction, $\frakM(A)$ is a special case of
$\frakA_k$.

(\mref{it:surj}) The assigned map is clearly well-defined and surjective.
It is an algebra homomorphism because the products in $\frakM(A)$ and
$\frakA_k$ are both defined by the Rota-Baxter relation~(\mref{eq:RB3}).
Alternatively, it can be proved by induction, as in item (\mref{it:ak}).

(\mref{it:eval}) follows from item (\mref{it:surj}) since a composition
of algebra homomorphisms is still an algebra homomorphism.
\end{proof}

\begin{coro}
If $F$ is an algebraic relation among elements $f_i, 1\leq i\leq n$
in $\frakM(A)_k$ for a given $k$,
then $\frakP_k(F)$ gives the same algebraic relation among the elements
$\frakP_k(f_i), 1\leq i\leq n$ in the MZV algebra $\frakA_k$ and
the same algebraic relation among elements $(\nu\circ \frakP_k)(f_i),
1\leq i\leq n$ of $\bfk$.
%In particular, any identity that
%holds for any element in a Rota--Baxter algebra gives an identity
%in a MZV\redt{$-$}algebra \redt{!!there was a $=$?!?}.
\mlabel{co:id}
\end{coro}

\subsection{Examples of MZV algebras}
As we have seen at the beginning of this section,
for the sum operator $P(f)(x)=\sum_{n=1}^\infty f(x+n)$,
we have
\begin{equation}
P_{(f_1,\cdots,f_k)}(x)=\sum_{n_1>n_2>\cdots >n_k\geq 1} f_1(x+n_1)f_2(x+n_2)
\cdots f_k(x+n_k)
\mlabel{eq:sumit}
\end{equation}
under suitable convergence condition.

\subsubsection{Multiple Hurwitz zeta functions and MZVs}

We let
$$A_H:=\{ f_s(x):=1/x^s \big | s\in \NN\}$$
with filtration given
by $s$. Then $A_H$ is a filtered subalgebra and generates a MZV algebra.
More precisely, we have
$$P(f_s)(x) =\sum_{n=1}^\infty \frac{1}{(x+n)^s}=\zeta(s;x+1)$$
where $\zeta(s;x)$ is the Hurwitz zeta function~\cite{Ca1}.
Evaluated at $x=0$, we obtain the Riemann zeta function
$\zeta(k)$.

We define a {\bf multiple Hurwitz zeta function} to be an
iteration
$$ \zeta(s_1,\cdots,s_k;x+1):=P_{(f_{s_k},\cdots, f_{s_1})}(x) =
\sum_{n_1>n_2>\cdots>n_k\geq 1} \frac{1}{(x+n_1)^{s_1}\cdots (x+n_k)^{s_k}} $$
with $s_1>1$.
So the MZV algebra $\frakA_H$ of level 2 generated by $A_H$ is the algebra of
multiple Hurwitz zeta functions.

Taking the evaluation map $\nu(f(x))=f(0)$,
we obtain the algebra of multiple zeta values
$$ \zeta(s_1,\cdots,s_k):=
\sum_{n_1>n_2>\cdots>n_k\geq 1} \frac{1}{n_1^{s_1}\cdots n_k^{s_k}}. $$
For example, for $P(f_2)(x)=\sum_{n=1}^\infty \frac{1}{(x+n)^2}=\zeta(2;x+1)$, by the
Rota-Baxter relation (\mref{eq:sum1}), we have
$P(f_2)P(f_2)=P(f_2P(f_2))+P(P(f_2)f_2)+P(f_2f_2)$. Since $f_2f_2=f_4$,
we obtain
$$ \zeta(2;x+1)\zeta(2;x+1)=2\zeta(2,2;x+1)+\zeta(4;x+1).$$
Evaluating at $x=0$, we obtain
\begin{equation}
\zeta(2)\zeta(2)=2\zeta(2,2)+\zeta(4).
\mlabel{eq:zetas}
\end{equation}
%Theorem~\mref{thm:MZVsubAlg}
%implies the well-known fact that the product of two MZVs
%of length $l_1$ resp. $l_2$ and depth $d_1$ resp. $d_2$ equals a
%linear combination of MZVs of length $l \leq l_1+l_2$ and depth
%$d=d_1+d_2$.
Further applications of Rota-Baxter algebras to MZVs will be given in the
next section.

\subsubsection{Multiple Lerch functions and MPLs}

Now let $s\in\NN$ and $z$ be a (complex) parameter with $|z|<1$.
Let
$$
    A_L:=\{ f_{s,z}(x,y):=z^y/x^s \big | s\in \NN, |z|<1\}.
$$
Then $A_L$ is a filtered algebra from the grading by $s$.
%\mfootnote{Assume $y$ integer for $z^y$ to make sense?}
We have
$$
    P(f_{s,z})(x,y):=\sum_{n=1}^\infty \frac{z^{y+n}}{(x+n)^s}
$$ and
$P(f_{s,z})(x,0)=\Phi(z,s,x+1)$ where $\Phi(z,s,x)$ is the Lerch
function~\cite{Ca1}. Evaluated at $x=0$, we obtain the
polylogarithm function
$$
    Li_s(z): = \sum_{n\geq 1} \frac{z^n}{n^s}.
$$
Further $Li_s(1)$ is the Riemann zeta function.

We define a {\bf multiple Lerch function} to be an iteration
$$
    \zeta(s_1,\cdots,s_k;x+1):=P_{(f_{s_1,z_1},\cdots, f_{s_k,z_k})}(x,y) =
    \sum_{n_1>n_2>\cdots>n_k\geq 1}
    \frac{z_1^{y+n_1}\cdots z_k^{y+n_k}}{(x+n_1)^{s_1}\cdots (x+n_k)^{s_k}}.
$$
So the MZV algebra $\frakA_L$ generated by $A_L$ is the algebra of
multiple Lerch functions.

When $y=x=0$, we obtain the multiple polylogarithms
$$
    Li_{s_1,\cdots,s_k}(z_1,\cdots,z_k):=
    \sum_{n_1>n_2>\cdots>n_k\geq 1} \frac{z_1^{n_1}\cdots z_k^{n_k}}{n_1^{s_1}\cdots n_k^{s_k}}.
$$

\subsubsection{Variations of $q$-MZVs}

There are several versions of $q$-MZVs and $q$-MPLs.
See~\cite{Zh,Brad} for $q$-multiple zeta values and polylogarithms.
See also~\cite{Zud} for another definition of $q$-MZVs.
They can all be defined as iterations of summation operators on certain functions.
%\mfootnote{More details can be provided.}

Fix $0<q<1$
Let $A_q$ be the subspace with basis
$$
    \left\{ {\frak q}_s(k):=\frac{q^{k(s-1)}}{[k]_q^s}\ \bigg |\ s\in \CC\right\}.
$$
Here $[k]_q=\frac{1-q^k}{1-q}$. Then $A_q$ is a filtered
subalgebra since
$$
    {\frak q}_s(k){\frak q}_t(k)={\frak q}_{s+t}(k)
    +(1-q){\frak q}_{s+t-1}(k).
$$
The MZV algebra $\frakA_{q,2}$ of level 2 generated by $A_q$, after
the evaluation map,  is the algebra of
$q$-MZVs studied by Zhao, Bradley, Kaneko, et. al.~\mcite{Zh,Brad,O-T}
consisting elements of the form
\begin{equation}
    \zeta_q(s_1,\cdots,s_d):=\sum_{k_1>\cdots > k_d> 0}
        \frac{q^{k_1(s_1-1)+\cdots+k_d(s_d-1)}} {[k_1]^{s_1}\cdots
            [k_d]^{s_d}} \mlabel{eq:qmzv}
\end{equation}
where $s_1>1$.

%{\bf Observation:} Suppose there is a relation that holds in every RBA
%and thus a relation among the multiple Hurwitz/Lerch functions. With
%suitable evaluations, we obtain a relation among MZVs or MPLs.
%\\
%{\bf Question:} How to find such relations?
%\\
%{\bf Answer (Rota, Cartier, 1970s):} Find relations (solution of
%word problems) in free RBAs.

\subsection{Rota--Baxter algebras and nested sums}
There is another {angle to} the Rota--Baxter algebra point of
view for MZVs. Let $A$ be a unitary ring and let $\NN$ be the set
of positive integers. Define $\cala:=\map(\NN, A)=A^\NN$ to be the
algebra of maps $f:\NN\to A$ with point-wise operations. Define a
linear operator
\begin{equation}
Z:=Z_\cala: \cala \to \cala, \quad
  Z[f](k):=\left \{ \begin{array}{ll} \disp{\sum_{i=1}^{k-1}}f(i), & k>1, \\
    0, & k=1. \end{array} \right .
\end{equation}
\begin{lemma}
$Z$ is a Rota--Baxter operator on $\cala$ of weight 1.
\mlabel{lem:sum}
\end{lemma}

Further we have, for $f_1,\cdots,f_n\in \cala$,
\begin{equation}
Z[f_1 Z[f_2 \cdots Z[f_n]\cdots ]](k+1)=
    \sum_{k\geq i_1>\cdots>i_n>0} f_1(i_1) \cdots f_n(i_n)
\end{equation}
and thus
\begin{equation}
\lim_{k\to \infty} Z[f_1 Z[f_2 \cdots Z[f_n]\cdots ]](k)=
    \sum_{i_1>\cdots>i_n>0} f_1(i_1) \cdots f_n(i_n)
\end{equation}
if the nested infinite sum on the right exists.

For example, the multiple zeta values are obtained as
$$
    \zeta(s_1,\cdots,s_n)=\sum_{i_1>\cdots>i_n>0} \frac{1}{i_1^{s_1} \cdots i_n^{s_n}}
    = \lim_{k\to \infty} Z\Big[\frac{1}{x^{s_1}}Z\Big[\frac{1}{x^{s_2}}\cdots
    Z\Big[\frac{1}{x^{s_n}}\Big]\cdots \Big]\Big](k).
$$

\subsection{Rota--Baxter algebras and iterated integrals}

We mentioned in the introduction the double-shuffle
relations for MZVs, corresponding to the sum and integral representations
of MZVs (see \mcite{I-K-Z} for more details). We have just seen that
the sum representation of MZVs is captured as a subquotient $A_H$
of a free Rota-Baxter algebra of weight 1 (Theorem~\mref{thm:MZVsubAlg}).
Similarly, the integral representation of MZVs is captured by a Rota-Baxter
algebra of weight zero, the integral operator considered in Eq (\mref{eq:int}).
For example the integral representation
$$\zeta(2)=\int_0^1 \frac{dx_1}{x_1}\int_0^{x_1}\frac{dx_2}{1-x_2}$$
is the evaluation at $x=1$ of
$$ I_{(1/x_1,1/(1-x_2))}:=\int_0^x \frac{dx_1}{x_1}\int_0^{x_1}\frac{dx_2}{1-x_2}$$
for the integration operator $I$ in Eq. (\mref{eq:int}) which is a Rota-Baxter
operator of weight zero. Since the product in a free Rota-Baxter algebra of
weight zero is given by the shuffle product, we have
\begin{eqnarray*}
 I_{(1/x_1,1/(1-x_2))}I_{(1/y_1,1/(1-y_2))}&=& I_{(1/x_1,1/y_1,1/(1-x_2),1/(1-y_2))}+
                            I_{(1/x_1,1/y_1,1/(1-y_2),1/(1-x_2))}\\
                         & &\quad +I_{(1/y_1,1/x_1,1/(1-x_2),1/(1-y_2))}
                    +I_{(1/y_1,1/x_1,1/(1-y_2),1/(1-x_2))}\\
                 & & \quad\; +I_{(1/x_1,1/(1-x_2),1/y_1,1/(1-y_2))}
                                +I_{(1/y_1,1/(1-y_2),1/x_1,1/(1-x_2))}\\
 &=& 4I_{(1/z_1,1/z_2,1/(1-z_3),1/(1-z_4))}+2I_{(1/z_1,1/(1-z_2),1/z_3,1/(1-z_4))}.
\end{eqnarray*}
Evaluated at $x=1$, we have
$$\zeta(2)\zeta(2)=4\zeta(3,1)+2\zeta(2,2).$$
Combining with Eq (\mref{eq:zetas}) we have
the famous relation
$$
  \zeta(3,1)=\frac{1}{4}\zeta(4).
$$
In general, this double-shuffle structure is reflected in the context
of Rota-Baxter algebras and needs to be further analyzed.
\medskip

{We will only
briefly mention the integral representation of q-MZVs, as they
appeared in~\mcite{Zh,Brad}. The $q$-analog of the Riemann
integral, named Jackson integral after its inventor Reverend
Jackson~\mcite{Ro5}, on a well chosen function algebra is given
by}
\begin{equation}
    J[f](x) := \int_{0}^{x} f(y) d_qy
            := (1-q)\:\sum_{n \ge 0} f(xq^n) xq^n.   \mlabel{JI}
\end{equation}
for $0<q<1$.
A key ingredient in the Jackson integral is the operator
\begin{equation}
    P_q[f](x):=\sum_{n>0} f(xq^n).         \mlabel{P}
\end{equation}
\begin{prop} \mcite{Ro4}
The maps $P_q$ and $\hat{P}_q:=id+P_q$ are Rota-Baxter
operators of weight $1$ and $-1$, respectively.
\end{prop}
It then follows that Jackson's integral satisfies the relation
\begin{equation}
    J[f]\: J[g] + (1-q)J[f \: g \: \mathrm{id}] = J \Big[J [f] \: g + f \: J[g] \Big], \mlabel{JB}
\end{equation}
where $\mathrm{id}$ is the identity map.

The q-analogs of MZVs of~\cite{Zh,Brad} have a
Jackson-integral representation.
% \allowdisplaybreaks{
%\begin{eqnarray}
%     \zeta_q(s_1, \cdots, s_m)=\!\! \underbrace{\hat{P}_q \: \bigg[\hat{P}_q \: \Big[\dots \hat{P}_q}_{s_1}\: \Big[
%                                                              \hat{\Omega}\big\{t;\prod_{j=1}^{1}q^{1-s_j}\big\}
%                                     \dots
%                                     \underbrace{\hat{P}_q \: \big[\hat{P}_q\:\big[ \dots \hat{P}_q\!\!
%                                                         \phantom{\bigg[}}_{s_m}\big[
%                                                              \hat{\Omega}\big\{t;\prod_{j=1}^{m}q^{1-s_j}\big\}\big]
%                                    \dots \big]\big]\Big] \dots \Big]\bigg](1).\nonumber \mlabel{PqMZV}
%\end{eqnarray}}
%where $\hat{\Omega}\{t;b\}:=\frac{t}{b-t}$.
Unfortunately, Theorem \mref{thm:MZVsubAlg} does not apply to the
integral representation of these $q$-MZVs due to the lack of a suitable
MZV algebra structure here. We plan to elaborate on this in a future work.

\section{Identities in Rota--Baxter algebras and MZVs}
\mlabel{sec:app} By Corollary~\mref{co:id}, once we have an
identity in free Rota--Baxter algebras, we can apply it to various
MZV algebras to get identities there. We give some examples of such
applications by using old and new results in Rota--Baxter
algebras.

\subsection{Spitzer's identity}
Let $(R,P)$ be a unitary commutative Rota--Baxter $\QQ$-algebra of
weight $1$. Then for $a \in R$, we have the Spitzer's identity~\mcite{Ro1}
\begin{equation}
 \exp\left (P(\log(1+at)) \right )
    =\sum_{i=0}^\infty t^i \underbrace{P\big( P( P( \cdots P}_{i\mbox{\rm -}{\rm times}}(a)a)\cdots  a) a\big)
    \mlabel{eq:si1}
\end{equation}
in the ring of power series $R[[t]]$.
We apply it to the MZV algebra of level 2 generated by $A_H$.
For $k\geq 2$, let $a=f_k:=1/x^k\in A_H$. With the summation operator
$P(f)(x)=\sum_{n\geq 1} f(x+n)$, we have
\begin{eqnarray*}
P(\log (1+at))&=& P\big(\sum_{i=1}^\infty \frac{(-1)^{i-1}}{i}
\left(\frac{t}{x^k} \right)^i\big)\\
&=& \sum_{i=1}^\infty \frac{(-1)^{i-1}}{i} P\left(\frac{1}{x^{ki}} \right) t^i\\
&=& \sum_{i=1}^\infty \frac{(-1)^{i-1}}{i} \zeta(ki;x+1) t^i
\end{eqnarray*}
Similarly, the right hand side of Eq (\mref{eq:si1}) becomes
$$ 1+\sum_{i=1}^\infty \zeta (\underbrace{k,\cdots, k}_{i};x+1)t^{i}.$$
So Evaluating at $x=0$, we get the well-known identity
$$\exp\left (\sum_{i=1}^\infty (-1)^{i-1} \zeta (ik)\frac{t^{i}}{i}\right )
= 1+\sum_{i=1}^\infty \zeta (\underbrace{k,\cdots,
k}_{i})t^{i}.$$
Applying Spitzer's identity to the subalgebra
$A_L$ gives similar relations of MPLs.

The non-commutative version of Spitzer's identity (\ref{eq:si1}) gives
Bogoliubov's formulae in perturbative renormalization of quantum
field theory~\mcite{E-G-K3}.

We now apply Spitzer's identity to the free Rota-Baxter
algebra $A=\sha(\CC[x])$ with its Rota-Baxter operator still denoted by
$P$. For $u,v\in A$, define
$$ u\star v= uP(v)+P(u)v+uv.$$
Then define
$$u^{\star\, n}=\underbrace{u\star\cdots \star u}_{n},\quad
\exp_\star(u)=\sum_{n\geq 0} u^{\star\, n}/n!.$$
We verify that $P(u)P(v)=P(u\star v)$ by definition of $P$
and hence $P(u)^n=P(u^{\star\, n})$, $n\geq 1$. Using this to rewrite
Spitzer's {identity}, the left hand side gives
$P(\exp_\star\log(1+xt))-P(1)+1.$
The right hand side gives
\begin{eqnarray*}
 1+P(xt)+P(xtP(xt))+\cdots
        & = & 1+1\ot xt + 1\ot xt\ot xt +\cdots \\
        & = & P(1 +xt+(xt)^{\ot 2}+\cdots)-P(1)+1.
\end{eqnarray*}Thus $\exp_*\log(1+xt)=1 +xt+(xt)^{\ot 2}+\cdots.$
This is a basic identity in~\mcite{I-K-Z}.

\subsection{Bohnenblust-Spitzer formula}
We quote the following theorem that first appeared in the same
paper~\mcite{Sp} where Spitzer published his above formula and
takes it present form in~\mcite{Ro3}.
\begin{theorem}
{\bf (Bohnenblust-Spitzer formula)} Let $A$ be a Rota--Baxter
algebra. Then
$$ \sum_{\sigma\in S_n} P (s_{\sigma(1)}P(s_{\sigma(2)}\cdots P(s_{\sigma(n)})\cdots ))
  =\sum_{\calt} (-1)^{n-|\calt|}
    \prod_{T\in\calt} (|T|-1)! P( \prod_{j\in T} s_j),\ n>0.$$
Here $s_j>1, 1\leq j\leq n$ and
$\calt$ runs through all unordered set partitions of $\{1,\cdots, n\}$.
\end{theorem}
Applying this identity to $A_H$, we obtain an identity of multiple Hurwitz
zeta functions. Then specialized at $x=0$, we have the following important
Partition Identity of Hoffman~\mcite{Ho2}.

\begin{coro} For $s_i>1$, $1\leq i\leq n$, we have
$$ \sum_{\sigma\in S_n} \zeta (s_{\sigma(1)},\cdots, s_{\sigma(n)})\\
    =\sum_{\calt} (-1)^{n-|\calt|}
    \prod_{T\in\calt} (|T|-1)!\zeta( \sum_{j\in T} s_j),\ n>0.$$
\end{coro}
For example, when $n=2$, we have the identity
$$P(s_1P(s_2))+P(s_2P(s_1))=-P(s_1s_2)+P(s_1)P(s_2),$$
in Rota--Baxter algebras, translated to
$$\zeta(s_1,s_2)+\zeta(s_2,s_1)=-\zeta(s_1+s_2)+\zeta(s_1)\zeta(s_2).$$
Applying Bohnenblust-Spitzer formula to $A_q$, we have
the $q$-analog of Hoffman's identity proved by Bradley~\mcite{Brad}.

\subsection{Congruences}
The following congruence relation is prove in~\mcite{Gu3}.
\begin{theorem} Let $p$ be a prime number.
For any $a_1\otimes \cdots \otimes a_n\in A^{\otimes n}$ in the free
Rota-Baxter algebra $\sha(A)$,
we have
$$
{\rm (Tensor/Graduate\ form\ of\ Freshman's\ Dream)} \quad
(a_1\otimes \cdots \otimes a_n)^p
    \equiv a_1^p\otimes \cdots \otimes a_n^p
    \mod p$$
in the sense that $p$ divides the coefficients of all other pure tensors
when the power on the left hand side is expressed as a linear combination
of pure tensors. Here the product on the left hand side is the product
in $\sha(A)$ defined in Eq (\mref{eq:shpr}).
%Let $(A,P)$ be a RBA and let $a_i\in A$, then
%$$P(a_1P(a_2\cdots P(a_n)\cdots ))^p \equiv P(a_1^pP(a_2^p \cdots P(a_n^p)\cdots ))
%$$ modulo $p$.
\end{theorem}
Compare with the well-known Freshman's Dream
$$ (x+y)^p\equiv x^p+y^p \mod p.$$
Applying the theorem to the MZV algebra $A_H$ of multiple Hurwitz series and
evaluating at $x=0$, we have
%\begin{proof}
%\begin{align*}
%(1\ot x)^p&=\sum_{i_1+\cdots+i_k=p} \binc{p}{i_1,\cdots,i_k}
%x^{i_1}\ot \cdots \ot x^{i_k}\\
%&\equiv 1\ot x^p \mod p.
%\end{align*}
%\end{proof}
\begin{coro}
$$ \zeta (s_1,\cdots,s_n)^p \equiv \zeta (ps_1,\cdots,ps_n) \mod p$$
in the sense that the coefficients of all other multiple zeta values on
the right hand side is a multiple of $p$.
\end{coro}

Similar congruences hold for MPLs.
\medskip

\noindent
{\bf Acknowledgements:} The second named author is supported in part by
NSF grant DMS-0505643 and a grant from Rutgers University Research
Council, and thanks IHES for hospitality.

%=========================================================================================
%=========================================================================================
%=========================================================================================
%\addcontentsline{toc}{section}{\numberline {}References}

\end{document}